\documentclass[11p]{article}
\usepackage[ansinew]{inputenc}
\usepackage{array}
\usepackage{color}
\usepackage{amsmath,amsthm,amscd,scalerel}
\usepackage{amsxtra}
\usepackage{amstext}
\usepackage{amssymb}
\usepackage{hyperref}
\usepackage{latexsym}
\usepackage{amsfonts,cite}
\usepackage{graphics,epstopdf}
\usepackage{epsfig}
\usepackage{amsfonts}
\usepackage{graphicx}
\usepackage{tcolorbox}
\usepackage{array,tabularx,colortbl}
\usepackage[framemethod=TikZ]{mdframed}

\providecommand{\U}[1]{\protect\rule{.1in}{.1in}}

\begin{document}
	
	\sloppy
	\newtheorem{thm}{Theorem}
	\newtheorem{cor}{Corollary}
	\newtheorem{lem}{Lemma}
	\newtheorem{prop}{Proposition}
	\newtheorem{eg}{Example}
	\newtheorem{defn}{Definition}
	\newtheorem{rem}{Remark}
	\newtheorem{note}{Note}
	\numberwithin{equation}{section}
	
	\thispagestyle{empty}
	\parindent=0mm
	
	\begin{center}
		{\large \textbf{Higher-order Hermite numbers: Properties and\\ applications to evolution problems}}\\ 					
		
		\vspace{0.20cm}
		
		{\bf Giuseppe Dattoli$^{1}$, \bf Subuhi Khan$^{2}$, Ujair Ahmad$^{3}$}\\
		\vspace{0.15cm}
		
		{$^{1}$ENEA - Frascati Research Center, Via Enrico Fermi 45, 00044 Frascati, Rome, Italy.}\\
		{$^{2,3}$Department of Mathematics, Aligarh Muslim University, Aligarh- 202001, India.}
			\footnote{$^{*}$This work has been done under Senior Research Fellowship (Ref No. 231610072319 dated:29/07/2023) awarded to the third author by University Grants Commission, New Delhi.}
		\footnote{$^{1}$pinodattoli@libero.it (G. Dattoli)}
		\footnote{$^{2}$subuhi2006@gmail.com (Subuhi Khan)}
		\footnote{$^{3}$ujairamu1998@gmail.com (Ujair Ahmad)(Corresponding Author)}
		
	\end{center}
	
	\begin{abstract}
		\noindent
			The operational calculus associated with Hermite numbers has been shown to be an effective tool for simplifying the study of special functions. Within this context, Hermite polynomials have been viewed as Newton binomials, with the consequent possibility of establishing previously unknown properties. In this article, this method is extended to study the lacunary Hermite polynomials and obtain novel results concerning their generating functions, recurrence relations, differential equations and certain integral transforms. The proposed method is systematically applied to a variety of evolution equations. Furthermore, this idea is extended to combinatorial interpretation of these polynomials, broadening their applicability in mathematical analysis and discrete structures.
	\end{abstract}
	\parindent=0mm
	\vspace{.25cm}
	
	
	\vspace{0.25cm}
	
	\section{Introduction}
	In 1861, John Blissard \cite{Blissard} introduced a novel notation for handling sums containing binomial coefficients. His approach involved expanding polynomials and then converting exponents into subscripts. For instance, the polynomial $(x+1)^n$ could be interpreted as $\sum_{k=0}^{\infty}\binom{n}{k}x_k$. Fifteen years later, Edouard Lucas \cite{Lucas,Lucas2} published his foundational papers on this subject, coining the terms ``symbolic method", or ``umbral calculus" \cite{Bell} . During the 20th century, Roman \cite{Roman}, Roman and Rota \cite{RomanRota}, Rota and Taylor \cite{RotaTaylor} referred to this systematically developed approach as ``classical umbral calculus", acknowledging Blissard's comprehensive work.\\
	
	This methodology has been widely adopted by researchers, often without detailed explanation \cite{Pascal}. While some applications have been relatively straightforward, others have involved more sophisticated treatments. Rigorous algebraic formulations of the method have been presented by Bell \cite{Bell3} and Temple \cite{Temple}, while Rota and Mullin \cite{RotaMulin} have characterized it in terms of linear operators.\\
	
Recent developments have demonstrated the remarkable versatility of umbral calculus, which continues to find new applications across diverse mathematical domains. This persistent significance has kept the subject as an active area of research interest, see, for example, \cite{DGLM,DGMR,RZAE,LD}. This article reveals an intermediate level that has been largely neglected in the literature, yet offers significant methodological advantages. The enhanced formulation is achieved through systematic absorption of Hermite numbers \cite{Wein,KK,DGLM}, which are important elements. These numbers serve as an intriguing and effective tool for exploring the theory of Hermite polynomials and other special polynomials from an operational perspective. The novelties brought by these methods become even more evident when they are applied to the study of higher-order Hermite polynomials \cite{HM,HM2,DLS}.\\ 
	
	
	
The study of Hermite polynomials is elegantly unified through the framework of umbral calculus. We begin with the two-variable Hermite polynomials $H_{n}(x,y)$, defined by the following series expansion \cite{Appell}:
\begin{equation}\label{hneq1}
	H_{n}(x,y) = n! \sum_{r=0}^{\left\lfloor\frac{n}{2}\right\rfloor} \frac{x^{n-2r} y^{r}}{r! (n-2r)!}.
\end{equation}
A significant property of $H_{n}(x,y)$ is its operational representation in the following form:
\begin{equation}\label{hneq68}
	H_{n}(x,y) =e^{y\partial_{x}^{2}}x^{n}.
\end{equation}
These polynomials generalize the classical Hermite polynomials $H_{n}(x)$, which are recovered via relation
	\begin{equation}\label{hneq61}
	H_{n}(x)=H_{n}(2x,-1).
\end{equation}
Relation \eqref{hneq61} may be used to verify the well known series representation of the classical Hermite polynomials $H_{n}(x)$:
\begin{equation}\label{hneq59}
	H_{n}(x) = n! \sum_{r=0}^{\left\lfloor\frac{n}{2}\right\rfloor} \frac{(-1)^{r} (2x)^{n-2r}}{r! (n-2r)!}.
\end{equation}
 We recall that Hermite numbers are the values of the Hermite polynomials evaluated at zero, that is
\begin{equation}\label{hneq60}
	H_{n}:=H_{n}(0),
\end{equation}
which admits the closed-form expression as \cite{Wein}:
\begin{equation}\label{hneq67}
	H_{n}=\dfrac{2^{n}\sqrt{\pi}}{\Gamma\left(\dfrac{1-n}{2}\right)},
\end{equation}
Expression \eqref{hneq67} generates the following sequence:
\begin{equation}\label{hneq55}
	H_{n}=\dfrac{(-1)^{n/2}n!}{\Gamma(n/2+1)}\left|\cos\dfrac{r\pi}{2}\right|=1, 0, -2, 0, 12, 0, -120,\cdots. \quad \text{[\url{(http://oeis.org/A067994)}]}
\end{equation}
To develop an umbral formulation, we define the umbral shift operator $\hat{h}_{\mu}^{r}$ as a function
\begin{equation}\label{hneq56}
	\hat{h}_{\mu}^{r}\;:\; \varphi(\mu)\to\varphi(\mu+r),
\end{equation}
where, the function 
\begin{equation}\label{hneq57}
	\varphi(\mu):=\varphi_{\mu}=\dfrac{\mu!}{\Gamma(\mu/2+1)}\left|\cos\dfrac{\mu\pi}{2}\right|
\end{equation}
is termed as umbral vacuum with initial state $\varphi(0):=\varphi_{0}=1$.\\

Remarkably, the two-variable Hermite polynomials $H_{n}(x,y)$ can be expressed as:
\begin{equation}\label{hneq58}
H_{n}(x,y)=(x+\sqrt{y}\hat{h}_{\mu})^{n}\varphi_{\mu}|_{\mu=0}.
\end{equation}
Expanding Newton binomial in \eqref{hneq58} and simplifying, the original series \eqref{hneq1} of $H_{n}(x,y)$ is obtained, which validates the umbral representation.\\

The formalism can be further streamlined by defining the umbral operator $\hat{h}$ acting on the initial state $\varphi_{0}$ such that
\begin{equation}\label{hneq4}
	\hat{h}^{r}\varphi_{0}=\varphi(r)=\varphi_{r}=\dfrac{r!}{\Gamma(r/2+1)}\left|\cos\dfrac{r\pi}{2}\right|,
\end{equation}
allowing us to rewrite umbral equation \eqref{hneq58} as
	\begin{equation}\label{hneq3}
		H_{n}(x,y)=(x+\sqrt{y}\hat{h})^{n}\varphi_{0}.
	\end{equation}

	As a consistency check, the classical Hermite polynomials $H_{n}(x)$ in this framework become 
	\begin{equation}\label{hneq63}
		H_{n}(x)=(2x+(-1)^{1/2}\hat{h})^{n}\varphi_{0}.
	\end{equation}
whose Newton binomial expansion yields series expansion \eqref{hneq59}.\\

Equation \eqref{hneq63} at $x=0$ provides following umbral representation of Hermite numbers $H_{n}$:
\begin{equation*}
	H_{n}:=H_{n}(0)=((-1)^{1/2}\hat{h})^{n}\varphi_{0},
\end{equation*}
which correctly generates sequence \eqref{hneq55}.\\

A natural extension is obtained by setting $x=0$ and $y=1$ in umbral equation \eqref{hneq3}, defining the second-order Hermite numbers, denoted by $_{2}h_{r}$:
\begin{equation}\label{hneq62}
	_{2}h_{r}:=H_{n}(0,1)={_{2}\hat{h}^{n}}\varphi_{0}={_{2}\varphi_{r}},
\end{equation}
corresponding to the sequence:
\begin{equation}\label{hneq2}
	_{2}h_{r}=\dfrac{r!}{\Gamma(r/2+1)}\left|\cos\dfrac{r\pi}{2}\right|=1, 0, 2, 0, 12, 0, 120,\cdots.
\end{equation}
To ensure clear distinction between Hermite numbers of different orders, we adopt the modified notation where the umbral operator for second-order Hermite numbers is denoted as $_2{\hat{h}}$ instead of $\hat{h}$, while its corresponding initial vacuum state is written as $_2{\varphi_{0}}$ rather than $\varphi_{0}$. This explicit notational convention will help maintain clarity when discussing operators and states associated with specific orders throughout this article.\\

	The operator $_{2}\hat{h}$ offers a powerful tool for computations. For instance, the Gaussian function is expressible as an ordinary exponential
	\begin{equation}\label{hneq5}
		e^{i\;{_{2}\hat{h}}x}{_2{\varphi}_{0}}=e^{-x^{2}}.
	\end{equation}

	The structure of this paper is organized as follows. In Section $2$, the utilities of the umbral approach are demonstrated through investigation of key properties and derivation of new results for second-order Hermite numbers $_2{h_{r}}$. In Section $3$, the framework is extended to establish analogous results for third-order Hermite numbers $_3{h_{r}}$ and their higher-order generalizations. Practical applications of these theoretical developments are presented in Section $4$. The article is concluded in Section $5$ with remarks on potential directions for future research.\\
	\section{Second order Hermite numbers}
	The use of standard calculus methods combined with the formalism outlined in introduction, leads to the derivation of infinite integrals. To illustrate the process, an infinite integral involving umbral form of Gaussian function is derived in the following example:
	\begin{eg} Let us evaluate the following infinite integral of Gaussian function:
		\begin{equation*}
			I(\nu)=\left|\int_{0}^{\infty}e^{-x^{2}}x^{\nu-1}dx\;{_2{\varphi}_{0}}\right|.
		\end{equation*}
	In view of equation \eqref{hneq5} and making the substitution $i\;_2{\hat{h}}x = -\zeta$ within the integral, it follows that
		\begin{equation*}
		I(\nu)=	\left|\int_{0}^{\infty}e^{i\;_2{\hat{h}}x}x^{\nu-1}dx\;{_2{\varphi}_{0}}\right|=\left|\left(-\dfrac{1}{i\;_2{\hat{h}}}\right)^{\nu}\int_{0}^{\infty}e^{-\zeta}\zeta^{\nu-1}d\zeta\;{_2{\varphi}_{0}}\right|,
		\end{equation*}
		
	which on using the integral definition of the gamma function simplifies to
		\begin{equation*}
			I(\nu)=\left|\left(-\dfrac{1}{i\;_2{\hat{h}}}\right)^{\nu}\Gamma(\nu)\;{_2{\varphi}_{0}}\right|.
		\end{equation*}
		
	Applying operator equation \eqref{hneq4} and simplifying, we have
		\begin{equation*}
			I(\nu)=\left|\dfrac{\Gamma(\nu)\Gamma(1-\nu)}{\Gamma\left(1-\dfrac{\nu}{2}\right)}\left|\cos\left(\dfrac{\pi\nu}{2}\right)\right|\right|.
		\end{equation*}
		Employing the standard reflection formula for gamma function,
		we find the intermediate result
		\begin{equation}\label{hneq64}
			I(\nu)=\left|\int_{0}^{\infty}e^{-x^{2}}x^{\nu-1}dx\;{_2{\varphi}_{0}}\right|=\dfrac{\pi\left|\cos\left(\dfrac{\pi\nu}{2}\right)\right|}{\Gamma\left(1-\dfrac{\nu}{2}\right)\left|\sin(\pi \nu)\right|},
		\end{equation}
		which on further simplification leads to
		\begin{equation}\label{hneq6}
		I(\nu)=	\left|\int_{0}^{\infty}e^{-x^{2}}x^{\nu-1}dx\;{_2{\varphi}_{0}}\right|=\dfrac{1}{2}\Gamma\left(\dfrac{\nu}{2}\right).
		\end{equation}
		
	
	\end{eg}
	We establish another infinite integral involving the operator $_2{\hat{h}}$ in the next example.
\begin{eg}
	  Consider the evaluation of the following integral:
	  \begin{equation}\label{hneq65}
	  	I=\left|\int_{-\infty}^{\infty}e^{i\;_2{\hat{h}}x^{2}}dx\;{_2{\varphi}_{0}}\right|
	  \end{equation}
	Using the standard Gaussian integral formula
	  	\begin{equation}\label{hneq54}
	  		\int_{-\infty}^{\infty}e^{-ax^{2}+bx}dx=\sqrt{\dfrac{\pi}{a}}e^{\dfrac{b^{2}}{4a}},
	  	\end{equation}
for $b=0$, we immediately obtain
	  	\begin{equation*}
	  		\left|\int_{-\infty}^{\infty}e^{i\;_2{\hat{h}}x^{2}}dx\;{_2{\varphi}_{0}}\right|=\sqrt{\pi}\;_2{\hat{h}}^{-1/2}\;{_2{\varphi}_{0}}.
	  	\end{equation*}
 	  	Application of operator action specified in \eqref{hneq4} leads to
	  	\begin{equation}\label{hneq8}
	  	I=	\left|\int_{-\infty}^{\infty}e^{i\;_2{\hat{h}}x^{2}}dx\;{_2{\varphi}_{0}}\right|=\sqrt{\pi}\dfrac{\Gamma\left(\dfrac{1}{2}\right)}{\Gamma\left(\dfrac{3}{4}\right)}\cos\left(\dfrac{\pi}{4}\right).
	  	\end{equation}
	  	which on employing the Legendre's duplication formula yields the following result:
	  	 \begin{equation}\label{hneq7}
	  	I=\left|\int_{-\infty}^{\infty}e^{i\;_2{\hat{h}}x^{2}}dx\;{_2{\varphi}_{0}}\right|=\dfrac{1}{2}\Gamma\left(\dfrac{1}{4}\right).
	  	\end{equation}
\end{eg}
The above example is immediately understood to be correct since
\begin{equation}\label{hneq9}
	e^{i\;_2{\hat{h}}x^{2}}{_2{\varphi}_{0}}=e^{-x^{4}}.
\end{equation}
Integrals in equations \eqref{hneq6} and \eqref{hneq7} are evaluated by treating the operator $_2{\hat{h}}$ as an ordinary algebraic quantity.\\

 The framework enables numerous results, including compact expressions for repeated derivatives of composite functions, as demonstrated in the following result:
\begin{prop}
	The Hermite numbers $_{2}h_{n}$ can be expressed in terms of the two-variable Hermite polynomials $H_{n}(x,y)$ as follows:
	\begin{equation}\label{hneq11}
		\sum_{s=0}^{\infty}\dfrac{_{2}h_{s+n}}{s!}x^{s}=H_{n}(2x,1)e^{x^{2}}.
	\end{equation} 
	\begin{proof}
		Considering the operator relation
		\begin{equation*}
			e^{x^{2}}=e^{_2{\hat{h}}\;x}{_2{\varphi}_{0}}.
		\end{equation*}
		Differentiating both sides of this equation $n$ times with respect to $x$ gives
		\begin{equation*}
			\partial_{x}^{n}e^{x^{2}}={_2{\hat{h}}}^{n}e^{_2{\hat{h}}\;x}{_2{\varphi}_{0}}.
		\end{equation*} 
		Expanding the exponential operator in the r.h.s. and applying operation \eqref{hneq62}, it follows that
		\begin{equation}\label{hneq10}
			\partial_{x}^{n}e^{x^{2}}=\sum_{s=0}^{\infty}\dfrac{_{2}h_{s+n}}{s!}x^{s}.
		\end{equation}
		The left-hand side of equation \eqref{hneq10} can be interpreted using the known result \cite{BDLSA}:
		\begin{equation}\label{hneq77}
			\partial_{x}^{n}e^{x^{2}}=H_{n}(2x,1)e^{x^{2}}.
		\end{equation}
	In view of representations \eqref{hneq10} and \eqref{hneq77}, assertion \eqref{hneq11} follows.
	\end{proof}
\end{prop}
In order to show that the error function $\text{erf}(x)$ can be expressed in terms of second-order Hermite numbers $_{2}h_{r}$, we prove the following result:
\begin{thm}
For the error function $\text{erf}(x)$, the following expansion in terms of the second-order Hermite numbers $_{2}h_{r}$ holds true:
	\begin{equation}\label{hneq12}
		\text{erf}(x)=\dfrac{2}{\sqrt{\pi}}\sum_{s=0}^{\infty}\dfrac{_{2}h_{s}\;i^{s}}{(s+1)!}x^{s+1}.
	\end{equation}
	\begin{proof}
	Since the error function $\text{erf}(x)$ is defined by integral
	\begin{equation*}
		\text{erf}(x)=\dfrac{2}{\sqrt{\pi}}\int_{0}^{x}e^{-\zeta^{2}}d\zeta.
	\end{equation*}
Therefore, using operator equation \eqref{hneq5} and integrating, we find
	\begin{equation*}
		\text{erf}(x)=\dfrac{2}{\sqrt{\pi}}(i\;_2{\hat{h}})^{-1}(e^{i\;_2{\hat{h}}\;x}-1)\;{_2{\varphi}_{0}},
	\end{equation*}
which on expanding the exponential and using operator \eqref{hneq62} yields expansion \eqref{hneq12}.
\end{proof}
\end{thm}
The umbral representation of Hermite numbers can be utilized more flexibly by examining the integral
\begin{equation}\label{hneq13}
	I(\alpha) = \int_{-\infty}^{\infty} e^{-x^{2}} e^{-\alpha\; _2{\hat{h}}^{1/2}x} dx \;{_2{\varphi}_{0}},
\end{equation}
whose integrand is illustrated in Figure 1 for different values of $\alpha$. By using the ordinary Gaussian integration rule \eqref{hneq54}, we obtain
\begin{equation}\label{hneq14}
	I(\alpha) = \sqrt{\pi} e^{_2{\hat{h}}\left(\frac{\alpha}{2}\right)^{2}}\; {_2{\varphi}_{0}} = \sqrt{\pi} e^{\left(\frac{\alpha}{2}\right)^{4}},
\end{equation}
which is significant for various reasons. For instance, replacing $\alpha$ by $2i\zeta$, integral \eqref{hneq14} can be viewed as the integral representation of the super-Gaussian function $e^{-x^{4}}$ of order $4$.
\begin{center}
		\includegraphics[width=10cm,scale=5]{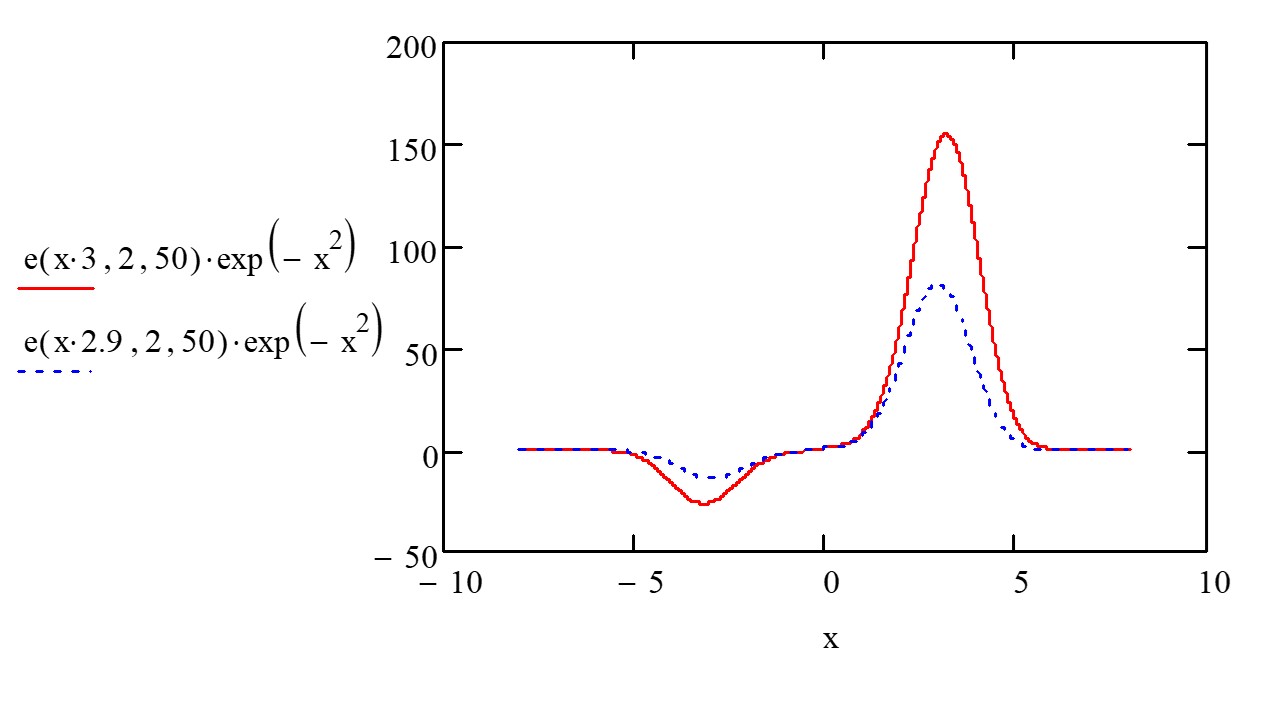}\\{\em \hspace{1cm}Figure 1 Graph of integrand in \eqref{hneq13} versus $x$ for $\alpha$ =$3$ (continuous line); $\alpha= 2.9$ (dotted line).}\\
\end{center}
A more computationally useful form of this solution can be obtained by the Gauss-Weierstrass transform \cite{BDLSA}, as shown in the following result:
\begin{prop}
For the	exponential operator $e^{y\partial_{x}^{4}}$, the following integral representation holds:
	\begin{equation}\label{hneq17}
		  e^{y\partial_{x}^{4}}=\dfrac{1}{\sqrt{\pi}}\int_{-\infty}^{\infty}e^{-\zeta^{2}}e^{-2\sqrt[4]{y}\partial_{x}\;{_2{\hat{h}}}^{1/2}\zeta}d\zeta\;{_2{\varphi}_{0}}.
	\end{equation}
	\begin{proof}
		Let us consider the expression
		\begin{equation*}
			F(x,y)=e^{y\partial_{x}^{4}}f(x),
		\end{equation*}
		which leads to an explicit solution for the evolution problem in terms of the following integral transform:
		\begin{equation}\label{hneq18}
			F(x,y)=\dfrac{1}{\sqrt{\pi}}\int_{-\infty}^{\infty}e^{-\zeta^{2}}e^{-2\sqrt[4]{y}\partial_{x}\;{_2{\hat{h}}}^{1/2}\zeta}f(\zeta)d\zeta\;{_2{\varphi}_{0}}.
		\end{equation}
		For $f(x)=1$, equation \eqref{hneq18} yields assertion \eqref{hneq17}.
	\end{proof}
\end{prop}
	\begin{rem}
		By making use of shift property of exponential operator $e^{\alpha\partial_{x}}$ in the integral \eqref{hneq18}, it follows that
		\begin{equation}\label{hneq19}
			F(x,y)=\dfrac{1}{\sqrt{\pi}}\int_{-\infty}^{\infty}e^{-\zeta^{2}}f(x-2\sqrt[4]{y}\;{_2{\hat{h}}}^{1/2}\zeta)d\zeta\;{_2{\varphi}_{0}},
		\end{equation}
		which provides numerically manageable form, assuming that the integral converges.
	\end{rem}

 The usefulness of this transform will become clearer in the context of higher-order Hermite polynomials.\\
 
The examples and results presented here illustrate the potency and reliability of the umbral formalism involving Hermite numbers and operators. In the forthcoming sections, we will discuss their role in more advanced topics. 
\section{Higher-order Hermite polynomials and associated numbers}
We recall that the third-order two-variable Hermite polynomials $H_{n}^{(3)}(x,y)$ are defined through the following series representation \cite{Appell,BDLSA}:
\begin{equation}\label{hneq22}
	H_{n}^{(3)}(x,y)=n!\sum_{r=0}^{\left\lfloor\frac{n}{3}\right\rfloor}\dfrac{x^{n-3r}y^{r}}{(n-3r)!r!}.
\end{equation}
The associated Hermite numbers are represented as:
\begin{align}\label{hneq23}
	_3{h}_r=\dfrac{r!}{\Gamma\left(\dfrac{r}{3}+1\right)}\left[2\left|\cos\left(\dfrac{r\pi}{3}\right)\right|-\left|\cos(r\pi)\right|\right]=1,0,0,6,0,0,360,0,0,60480,1.996.10^{7},\cdots,
\end{align}
and the relevant Newton binomial form is expressed as:
\begin{equation}\label{hneq24}
	H_{n}^{(3)}(x,y)=(x+\sqrt[3]{y}	\;{_3{\hat{h}}})^{n}\;	_3{\varphi}_0.
\end{equation}
where $_3{\hat{h}}$ is an umbral operator acting on $_3{\varphi}_{0}$ such that
\begin{equation}\label{hneq26}
	_3{\hat{h}}^{r}{_3{\varphi}_{0}}={_3{h}_{r}}.
\end{equation}
It is accordingly evident that the formalism, we have developed in the preceding section, holds almost unchanged for the forthcoming discussion. While the framework is conceptually straightforward, some notable implications are worth highlighting.\\
 
The third-order three-variable Hermite polynomials $H_{n}^{(3)}(x,y,z)$ are defined by the following series and generating function \cite{DAT}:
 \begin{equation}\label{hneq30}
 	H_{n}^{(3)}(x,y,z)=n!\sum_{r=0}^{\left\lfloor\frac{n}{3}\right\rfloor}\dfrac{H_{n-3r}(x,y)z^{r}}{(n-3r)!r!},
 \end{equation}
 and
 \begin{equation}\label{hneq70}
 	\sum_{r=0}^{\infty}\dfrac{t^{r}}{r!}H_{r}^{(3)}(x,y,z)=e^{xt+yt^{2}+zt^{3}}
 \end{equation}
 respectively.
 \begin{thm}
 For the third-order three-variable Hermite polynomials $H_{n}^{(3)}(x,y,z)$, the following series expansion holds true:
 \begin{equation}\label{hneq29}
 	H_{n}^{(3)}(-3x^{2}, -3x, -1)e^{-x^{3}}=(-1)^{n}\sum_{s=0}^{\infty}\dfrac{{_3{h}_{s+n}}}{s!}(-x)^{s},
 \end{equation}
 \begin{proof}
 	Let us consider
 	\begin{equation}\label{hneq25}
 		e^{-x^{3}}=e^{-{_3{\hat{h}}}\;x}{_3{\varphi}_{0}}.
 	\end{equation}
 	Differentiating above equation repeatedly $n$ times with respect to $x$, we obtain
 	\begin{equation*}
 		\partial_{x}^{n}e^{-x^{3}}=(-1)^{n}{	_3{\hat{h}}^{n}}\;e^{-{_3{\hat{h}}}\;x}{_3{\varphi}_{0}}.
 	\end{equation*}
 	Expanding the exponential term and applying operator action \eqref{hneq26}, we find
 	\begin{equation}\label{hneq27}
 		\partial_{x}^{n}e^{-x^{3}}=(-1)^{n}\sum_{s=0}^{\infty}\dfrac{{_3{h}_{s+n}}}{s!}(-x)^{s}.
 	\end{equation}
 	Since, we have the following identity \cite{BDLSA}:
 	\begin{equation}\label{hneq28}
 		\partial_{x}^{n}e^{-x^{3}}=H_{n}^{(3)}(-3x^{2}, -3x, -1)e^{-x^{3}},
 	\end{equation}
therefore, equating the $nth$ derivatives from equations \eqref{hneq27} and \eqref{hneq28}, assertion \eqref{hneq29} is obtained.
 \end{proof}
 \end{thm}
 
The utility of this technique is evident; all the properties of higher-order Hermite polynomials are reduced to those of ordinary Newton binomials, which makes the interplay between elementary and more advanced forms of special polynomials (or functions) extremely straightforward.\\

In the foregoing analysis, we utilized circular functions to introduce Hermite numbers for $m=2$ and $m=3$. Their use becomes rather plethoric for higher orders, as noted in \cite{DLS}. To address this, we adopt the following notation:
\begin{equation}\label{hneq38}
	_{m}h_{r}=\dfrac{r!}{\Gamma\left(\dfrac{r}{m}+1\right)}\delta_{\left\lceil \textstyle\frac{r}{m} \right\rceil, r}.
\end{equation}
For $m=4$, we have
\begin{equation*}
_{4}h_{r}	=1,0,0,0,24,0,0,0,2.016.10^{4},0,0,0,7.983.10^{7},\cdots.
\end{equation*}
The associated umbral operator is defined analogously as
\begin{equation}\label{hneq39}
	_m{\hat{h}}^{r}\;_m{\varphi_{0}}={	_{m}h_{r}}.
\end{equation}
The $m$-th order Hermite polynomials are given by
\begin{equation}\label{hneq40}
	H_{n}^{m}(x,y)=n!\sum_{r=0}^{\left\lfloor\frac{n}{m}\right\rfloor}\dfrac{x^{n-mr}y^{r}}{r!\;(n-mr)!}=(x+\sqrt[m]{y}\;_m{\hat{h}})^{n}{_m{\varphi_{0}}},
\end{equation}
while the non-lacunary version \cite{DAT} reads
\begin{equation}\label{hneq41}
	H_{n}^{m,m-1,\cdots,2}(x_{1},x_{2},\cdots,x_{m})=\left(x_{1}+\sum_{s=2}^{m}\sqrt[s]{x_{s}}\;_m{\hat{h}}\right)^{n}\prod_{s=2}^{m}\;_s{\varphi_{0}},
\end{equation}
which is the Newton multi-nomial form of the $m$-variable Hermite polynomials. This can be further generalized to Bell polynomials \cite{Bell}; see, also (\cite{Rioden}; pp. 35–38, 49, 142).\\

Exponential operators involving $m$-th order derivatives are consequently reduced to shift operators as follows:
\begin{equation}\label{hneq42}
	e^{y\partial_{x}^{m}}=e^{\sqrt[m]{y}\;_m{\hat{h}}\;\partial_{x}}{_m{\varphi_{0}}},
\end{equation}
thereby simplifying algebraic computations with higher-order derivative exponential operators.\\

To illustrate, let us consider the following result:

\begin{thm}
	For the triple lacunary Hermite polynomials $H_{3r}(x,y)$, the following generating relation holds:
	\begin{equation}\label{hneq47}
		\sum_{r=0}^{\infty}\dfrac{t^{r}}{r!}H_{3r}(x,y)=\;_3{E(x,y,t)}e^{tx^{3}},
	\end{equation}
	where
		\begin{equation}\label{hneq49}
		_3{E(x,y,t)}=\sum_{r=0}^{\infty}H_{2r}^{(3)}(3x^{2}t,3xt^{2},t)\dfrac{y^{r}}{r!}.
	\end{equation}
	\begin{proof}
		In view of operational definition \eqref{hneq68} of $H_{n}(x,y)$, the triple lacunary Hermite generating function can be written in the following form:
		\begin{equation}\label{hneq69}
			\sum_{r=0}^{\infty}\dfrac{t^{r}}{r!}H_{3r}(x,y)=e^{y\partial_{x}^{2}}e^{x^{3}t}.
		\end{equation}
Now, applying identity \eqref{hneq42} for $m=2$, the above generating function is expressed as:
		\begin{equation}\label{hneq43}
			\sum_{r=0}^{\infty}\dfrac{t^{r}}{r!}H_{3r}(x,y)=e^{\sqrt{y}\;_2{\hat{h}}\;\partial_{x}}e^{tx^{3}}{_2{\varphi_{0}}}.
		\end{equation}
		Application of shift exponential operator in the right hand side of equation \eqref{hneq43} gives
		\begin{equation}\label{hneq45}
			\sum_{r=0}^{\infty}\dfrac{t^{r}}{r!}H_{3r}(x,y)=e^{t(x+\sqrt{y}\;_2{\hat{h}})^{3}}\;_2{\varphi_{0}},
		\end{equation}
		which can be further rewritten as:
		\begin{equation}\label{hneq46}
			\sum_{r=0}^{\infty}\dfrac{t^{r}}{r!}H_{3r}(x,y)=e^{tx^{3}}e^{\left(3\sqrt{y}\;_2{\hat{h}}\;x^{2}+3y\;_2{\hat{h}}^{2}\;x+y^{3/2}\;_2{\hat{h}}^{3}\right)t}\;_2{\varphi_{0}}.
		\end{equation}
		 Now, define the the polynomials
		 \begin{equation}\label{hneq66}
		 	_3{E(x,y,t)}=e^{\left(3\sqrt{y}\;_2{\hat{h}}\;x^{2}+3y\;_2{\hat{h}}^{2}\;x+y^{3/2}\;_2{\hat{h}}^{3}\right)t}\;_2{\varphi_{0}}.
		 \end{equation}
		 Using generating function \eqref{hneq70} for three-variable Hermite polynomials, we define
		 \begin{equation}\label{hneq48}
		 	_3{E(x,y,t)}=e^{\left(3\sqrt{y}\;_2{\hat{h}}\;x^{2}+3y\;_2{\hat{h}}^{2}\;x+y^{3/2}\;_2{\hat{h}}^{3}\right)t}\;_2{\varphi_{0}}=\sum_{r=0}^{\infty}\dfrac{y^{r/2}\;_2{\hat{h}}^{r}}{r!}H_{r}^{(3)}(3x^{2}t,3xt^{2},t),
		 \end{equation}
		 which on account of the properties of the umbral operator $_2{\hat{h}}$ yields assertion \eqref{hneq47}, where $_3{E(x,y,t)}$ is defined in equation \eqref{hneq49}.
	\end{proof}
\end{thm}
\begin{note}
	The generating function \eqref{hneq43} is further written as:
	\begin{equation*}
		\sum_{r=0}^{\infty}\dfrac{t^{r}}{r!}H_{3r}(x,y)=e^{\sqrt{y}\;_2{\hat{h}}\;\partial_{x}}e^{\sqrt[3]{t}\;_3{\hat{h}}\;x}\;_3{\varphi_{0}}\;_2 {\varphi_{0}},
	\end{equation*}
		which on using the shift operator property, simplifies to
	\begin{equation}\label{hneq44}
		\sum_{r=0}^{\infty}\dfrac{t^{r}}{r!}H_{3r}(x,y)=e^{\sqrt[3]{t}\;_3{\hat{h}}(x+\sqrt{y}\;_2{\hat{h}})}\;_3{\varphi_{0}}\;_2{\varphi_{0}}.
	\end{equation}
\end{note}

The range of convergence of the series on the right-hand side of equation \eqref{hneq49} is rather narrow and indeed
$|t,y|<0.2$.\\

 An example of the behaviour of the function in \eqref{hneq47} with respect to $x$ is given in Figure $2$.
\begin{center}
	\includegraphics[width=10cm,scale=5]{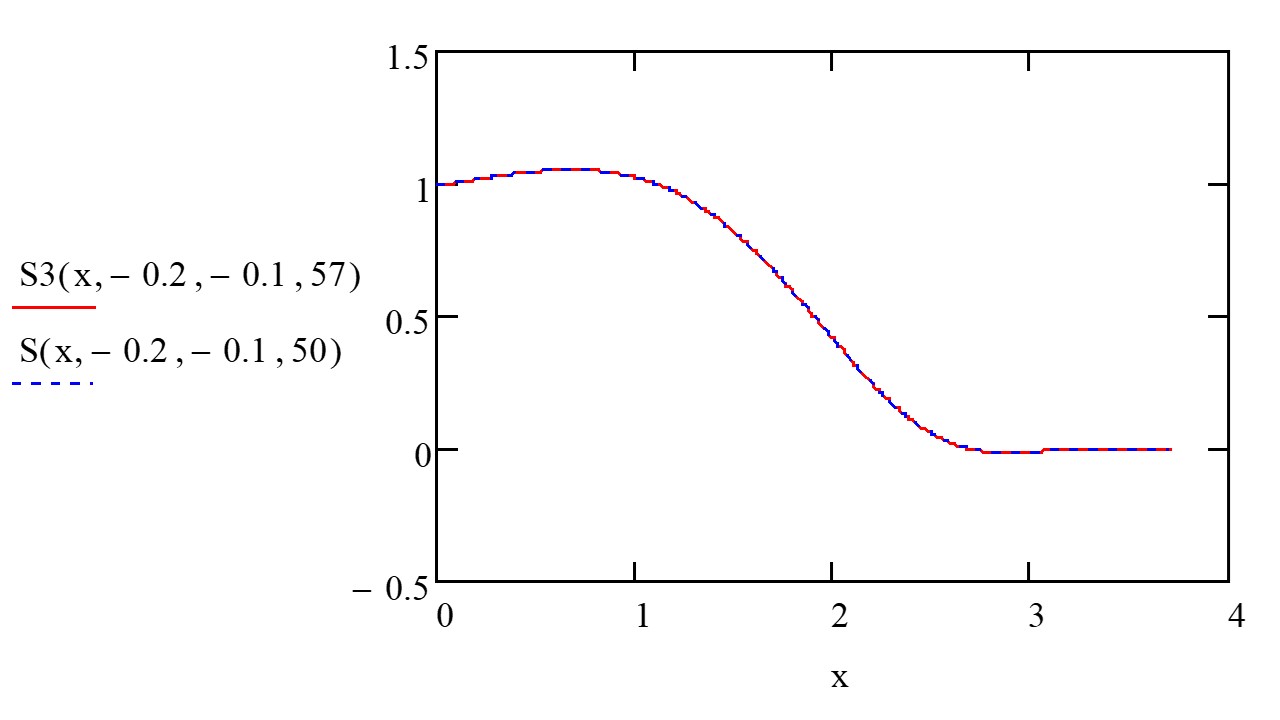}\\{\em \hspace{1cm}Figure 2 Comparison between equations \eqref{hneq48} and \eqref{hneq49} versus $x$ for $y=-0.2$, $t=-0.1$}\\
\end{center}
It is furthermore worth stressing that a result similar to \eqref{hneq47} has been obtained in \cite{GJ} using an umbral formalism analogous to ours but limited to second-order Hermite numbers.\\

Preceding sections in the article have established key findings and methodological developments that offer both fundamental insights and practical implications. These results not only advance theoretical understanding but also reveal promising opportunities for practical applications, which will be explored in the following section.
\section{Applications}
 In this section, the implementation of the established outcomes within certain domains is investigated. The translation from conceptual validation to practical application emphasizes the broader relevance and transformative potential of the results and developments.\\

An application of result \eqref{hneq14} is the solution of the fourth order heat equation \cite{HM2}:
\begin{equation}\label{hneq15}
	\partial_{y}F(x,y)=\partial_{x}^{4}F(x,y),
\end{equation}
with the initial condition
\begin{equation*}
	F(x,0)=g(x).
\end{equation*}
The solution of the above Cauchy problem can be formally expressed as:
\begin{equation}\label{hneq16}
	F(x,y)=e^{y\partial_{x}^{4}}F(x).
\end{equation}
Another area of exploration involves the Fourier-Gabor transform \cite{Gabor,JJD}, defined as:
\begin{equation}\label{hneq20}
	G(\tau,\omega) = \int_{-\infty}^{\infty} x(t) e^{-\pi (t-\tau)^{2} - i\omega t} dt.
\end{equation}
When expanded using the umbral formalism in terms of the two-variable Hermite polynomials $H_{n}(x,y)$, this transform yields the following result:
\begin{thm}
	The following representation of Fourier-Gabor transform in terms of two-variable Hermite polynomials holds:
	\begin{equation}\label{hneq21}
		G(\tau,\omega) = e^{-i\omega \tau} \sum_{n=0}^{\infty} A_{n}(\tau,\omega),
	\end{equation}
	where
	\begin{equation*}
		A_{n}(\tau,\omega) := (-i)^{n} \frac{H_{n}(\omega,\pi)}{n!} \int_{-\infty}^{\infty} x(t) (t - \tau)^{n} dt.
	\end{equation*}
	\begin{proof}
		Making use of umbral form of Gaussian function from equation \eqref{hneq5} into the right-hand integral of Fourier-Gabor transform \eqref{hneq20}, it follows that
		\begin{equation*}
			G(\tau,\omega) = \int_{-\infty}^{\infty} x(t) e^{i\;_2{\hat{h}}\sqrt{\pi} (t - \tau) - i\omega t} dt \varphi_{0},
		\end{equation*}
		which can be expressed as:
		\begin{equation*}
			G(\tau,\omega) = e^{-i\omega \tau} \int_{-\infty}^{\infty} x(t) e^{i(_2{\hat{h}}\sqrt{\pi} - \omega)(t - \tau)} dt \varphi_{0}.
		\end{equation*}
		Expanding the exponential term in the above integral, assertion \eqref{hneq21} is proved.
	\end{proof}
\end{thm}
Further important aspect of the formalism is its application to evolution problems. We note indeed that the polynomials in equation \eqref{hneq22} satisfy the third-order heat equation
\begin{equation}\label{hneq31}
	\partial_{y}F(x,y)=\partial_{x}^{3}(x,y),
\end{equation}
with the initial condition
\begin{equation*}
	F(x,0)=x^{n}.
\end{equation*}
The corresponding evolution operator can be expressed as:
\begin{equation}\label{hneq32}
	e^{y\partial_{x}^{3}}=e^{\sqrt[3]{y}\;{_3{\hat{h}}}\partial_{x}}{_3{\varphi}_{0}}.
\end{equation}
Exponential operators involving third-order derivatives are usually treated using the Airy transform method \cite{VS,ZK}. In this context, the reduction of the derivative order in the exponential is achieved through the identity \cite{Widder}
\begin{equation}\label{hneq33}
	e^{\lambda x^{3}}=\int_{-\infty}^{\infty}\exp\left\{\sqrt[3]{3y}\;xt\right\}A_{i}(t)dt\;;\quad\quad\lambda,~\Re{(x)}>0,
\end{equation}
where $A_{i}(t)$ denotes the Airy function \cite{VS}, defined as:
\begin{equation}\label{hneq34}
	A_{i}(t)=\dfrac{1}{2\pi}\int_{-\infty}^{\infty}\exp\left\{\dfrac{\iota}{3}\zeta^{3}+\iota t\zeta\right\}d\zeta.
\end{equation}
Replacing $x$ with $\partial_{x}$ in equation \eqref{hneq33}, we arrive at an alternate form of the evolution operator given in equation \eqref{hneq32}.\\

As previously noted, the solution of higher-order diffusive equations 
\begin{equation}\label{hneq50}
	\partial_{y}F(x,y)=\partial_{x}^{m}F(x,y),
\end{equation}
with the initial condition
\begin{equation*}
	F(x,0)=g(x)
\end{equation*}
can be obtained by taking advantage of the formalism, we have outlined. Therefore, let us take
\begin{equation}\label{hneq51}
	F(x,y)=e^{y\partial_{x}^{m}}g(x)=e^{\sqrt[m]{y}\;_m{\hat{h}}\;\partial_{x}}g(x)\;_m{\varphi_{0}}.
\end{equation}
Applying the Fourier transform property
\begin{align}\label{hneq52}
	g(x)=\dfrac{1}{\sqrt{2\pi}}\int_{-\infty}^{\infty}\tilde{g}(k)e^{ikx}dk;
\end{align}
\begin{align}\label{hneq53} \tilde{g}(k)=\dfrac{1}{\sqrt{2\pi}}\int_{-\infty}^{\infty}g(x)e^{-ikx}dk,
\end{align}
we obtain
\begin{equation*}
	F(x,y)=e^{\sqrt[m]{y}\;_m{\hat{h}}\;\partial_{x}}g(x)\;_m{\varphi_{0}}=\dfrac{1}{\sqrt{2\pi}}\int_{-\infty}^{\infty}\tilde{g}(k)e^{\sqrt[m]{y}\;_m{\hat{h}}\;\partial_{x}}e^{ikx}dk\;_{m}{\varphi_{0}},
\end{equation*}
Further using the shift property of the exponential operator, it follows that
\begin{equation*}
	F(x,y)=\dfrac{1}{\sqrt{2\pi}}\int_{-\infty}^{\infty}\tilde{g}(k)e^{ik(x+\sqrt[m]{y}\;_m{\hat{h}})}dk\;_{m}{\varphi_{0}}.
\end{equation*}
Examples of solutions for $m=2,4$ given in Figure $(3a)$, demonstrate their agreement with convention numerical methods. Additionally,  for $m=3$ (Figure $(3b)$), the solutions coincide with those reported in \cite{BDS2}.
\begin{center}
	\includegraphics[width=10cm,scale=5]{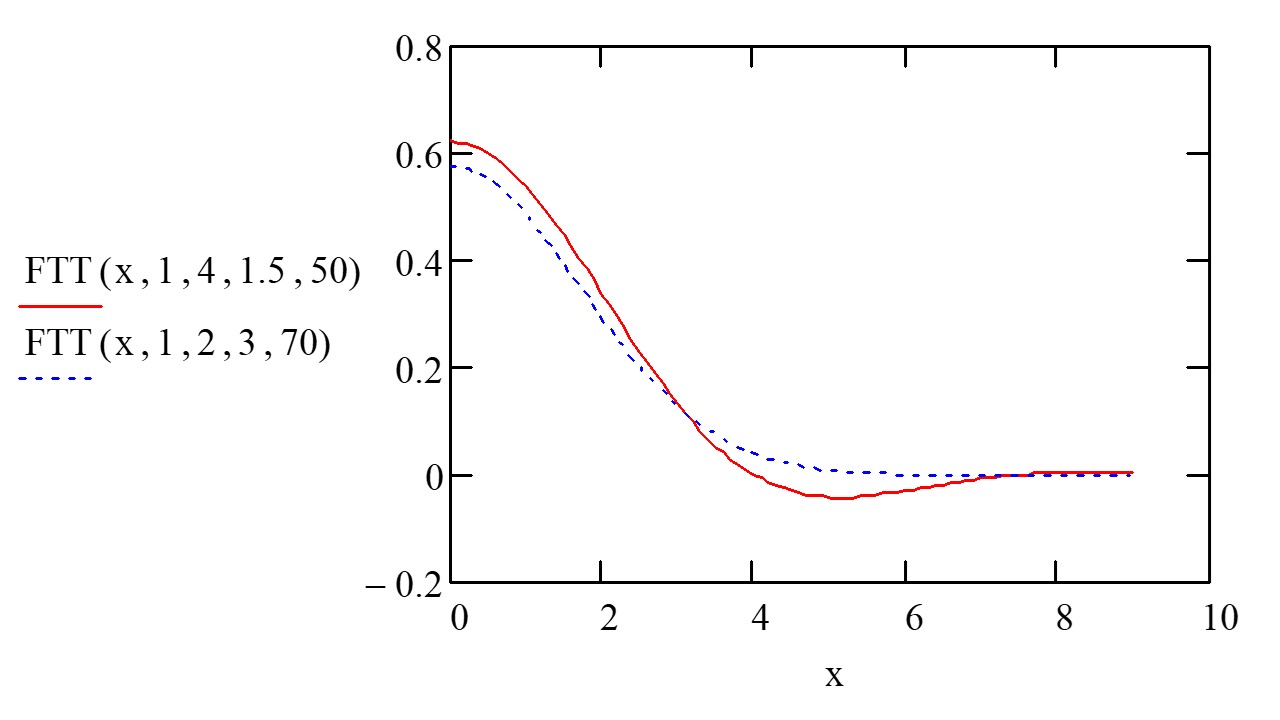}\\{\em \hspace{1cm}Figure 3(a) Solution of equation \eqref{hneq50} for $m=4,$ $y=1$ (continuous line); $m=2,$ $y=1$ (dot line)
	}\\
\end{center}
\begin{center}
	\includegraphics[width=10cm,scale=5]{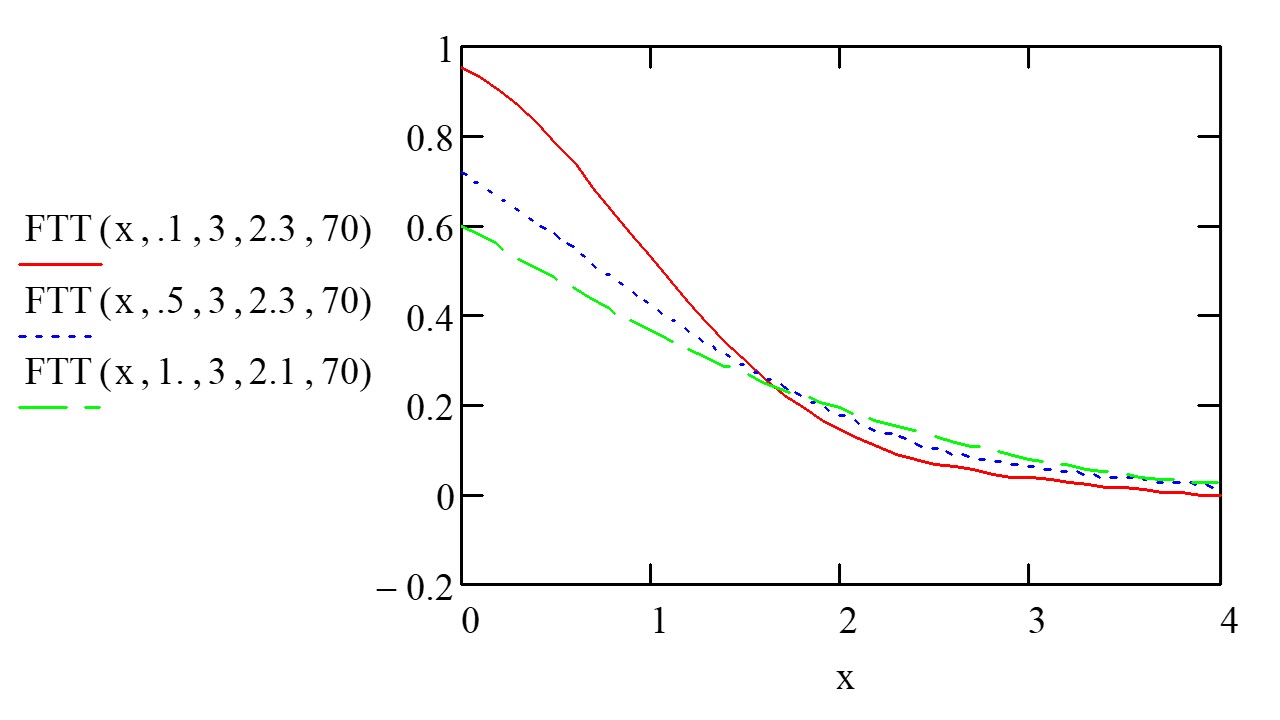}\\{\em \hspace{1cm}Figure 3(b) Solution of equation \eqref{hneq50} for $m=3,$ $y=0.1$ (continuous line); $m=3, y=0.5$ (dot line); $m=3,y=1$ (dash line)
	}\\
\end{center}
In the concluding section, we summarize the key outcomes and propose promising directions for future research.
\section{Concluding remarks}
This article explores the theory of Hermite polynomials, extending beyond the conventional two-variable case to investigate higher-order and multi-dimensional Hermite polynomials and numbers, as well as their relationship with generalized heat equations. The study presents a novel perspective by interpreting Hermite polynomials in terms of Newton binomials. A rigorous definition of umbral operators is provided, contextualized within the framework of the generalized Borel transform \cite{DPS}.\\

 It introduces higher-order Hermite numbers and establishes a previously unexplored connection between the Airy transform and Hermite polynomials. Additionally, the work demonstrates how higher-order heat equations can be systematically reduced to first-order partial differential equations, thereby offering a unified framework that enhances the applicability of integral transforms in solving such problems. The proposed methodology is systematically applied to a variety of evolution equations, demonstrating its utility in practical applications.\\
 
Before concluding this article, we emphasize that the three-variable Hermite polynomials \eqref{hneq30} can also be expressed as a Newton multi-nomial \cite{ENCYC}:
\begin{equation}\label{hneq35}
	H_{n}^{(3)}(x,y,z)=n!\sum_{r=0}^{\left\lceil\dfrac{n}{3}\right\rceil}\dfrac{(x+\sqrt{y}\;{_2\hat{h}})^{n-3r}z^{r}}{(n-3r)!r!}{_2{\varphi}_{0}},
\end{equation}
which can be further rewritten as:
\begin{equation}\label{hneq37}
	H_{n}^{(3)}(x,y,z)=\left(x+\sqrt{y}\;{_2\hat{h}}+\sqrt[3]{z}\;{_3\hat{h}}\right)^{n}{_2{\varphi}_{0}}{_3{\varphi}_{0}}.
\end{equation}
Further, equation \eqref{hneq37} gives
\begin{equation}\label{hneq36}
	H_{n}^{(3)}(x,y,z)=\sum_{k_{1}+k_{2}+k_{3}=n}\binom{n}{k_{1},k_{2},k_{3}}x^{k_{1}}\left(\sqrt{y}\;{_2\hat{h}}\right)^{k_{2}}\left(\sqrt[3]{z}\;{_3\hat{h}}\right)^{k_{3}}{_2{\varphi}_{0}}{_3{\varphi}_{0}},
\end{equation}
where
\begin{equation*}
	\binom{n}{k_{1},k_{2},k_{3}}=\dfrac{n!}{k_{1}!\;k_{2}!\;k_{3}!}.
\end{equation*}
The above results open the possibility of combinatorial interpretation (\cite{Stan} Vol. ~1, 2nd~ Ed., ~\S1.2) of the multi-variable Hermite polynomials.\\

  We have also outlined the usage of Hermite numbers in umbral form and studied extension of the method to super-Gaussian and to new forms of integral transforms. This last aspect will be developed in a forthcoming article to the study of super-Gaussian beams in optics \cite{BLMT,Gori}.\\

  \textbf{Conflict of Interest: }\\
  The authors declare that they have no conflict of interest.\\

  \textbf{Data availability:}\\
  Not applicable.
	
\end{document}